\documentclass[12pt]{amsart}
\usepackage{amssymb}
\usepackage{amscd}
\usepackage{a4}
\newcommand{\A}{\mathcal{A}}
\newcommand{\AG}{\A*_\alpha G}
\newcommand{\AGG}{\left(\AG\right)\# k[G]^*}
\newcommand{\B}{\mathcal{B}}
\newcommand{\AH}{\underline{\A \# H}}

\newcommand{\NN}{\mathbb{N}}
\newcommand{\rhu}{\rightharpoonup}
\newcommand{\lhu}{\leftharpoonup}
\newcommand{\End}[1]{\mathrm{End}_k({#1})}
\newcommand{\Ker}[1]{\mathrm{Ker}({#1})}
\renewcommand{\Im}[1]{\mathrm{Im}({#1})}
\newcommand{\e}{\mathbf{e}}
\renewcommand{\o}{\overline}
\newtheorem{thm}{Theorem}[section]
\newtheorem{prop}[thm]{Proposition}
\newtheorem{lem}[thm]{Lemma}
\newtheorem{cor}[thm]{Corollary}
\title{Duality for partial group actions}
\author{Christian Lomp}
\address{University of Porto, Department of pure Mathematics, Rua Campo Alegre 687, 4169-007 Porto (PORTUGAL)}
\email{clomp@fc.up.pt}

\thanks{This work was partially supported by Centro de Matem\'{a}tica da Universidade do Porto (CMUP), financed by FCT
(Portugal) through the programs POCTI (Programa Operacional Ci\^{e}ncia, Tecnologia, Inova\c{c}\~{a}o) and POSI
(Programa Operacional Sociedade da Informa\c{c}\~{a}o), with national and European community structural funds.}

\begin{document}

\maketitle

\begin{abstract}
Given a finite group $G$ acting as automorphisms on a ring $\A$, the skew group ring $\A*G$ is an important tool for studying the structure of $G$-stable ideals of $\A$. The ring $\A*G$ is $G$-graded, i.e.$G$ coacts on $\A*G$. The Cohen-Montgomery duality says that the smash product $\A*G \# k[G]^*$ of $\A*G$ with the dual group ring $k[G]^*$ is isomorphic to the full matrix ring $M_n(\A)$ over $\A$, where $n$ is the order of $G$. In this note we show how much of the Cohen-Montgomery duality carries over to partial group actions in the sense of R.Exel. In particular we show that the smash product $\AGG$ of the partial skew group ring $\AG$ and $k[G]^*$ is isomorphic to a direct product of the form $K\times \e M_n(\A) \e$  where $\e$ is a certain idempotent of $M_n(\A)$ and $K$ is a subalgebra of $\AGG$. Moreover $\AG$ is shown to be isomorphic to a separable subalgebra of $\e M_n(\A) \e$. We also look at duality for infinite partial group actions and for partial Hopf actions.
\end{abstract}

\section{Introduction}

Let $k$ be a commutative unital ring and $\A$ a untial $k$-algebra. Given a finite group $G$ acting as $k$-linear automorphisms on $\A$, Cohen and Montgomery showed in \cite{CohenMontgomery} that the smash product $\A*G \# k[G]^*$ of the skew group ring $\A*G$ and the dual group ring $k[G]^*=\mathrm{Hom}(k[G],k)$ is isomorphic to the full matrix ring $M_n(\A)$ over $\A$, where $n$ is the order of $G$. 

R.Exel introduced in \cite{exel} the notion of a partial group action on a $k$-algebra: $G$ acts partially on $\A$ by a family $\{\alpha_g :D_{g^{-1}} \rightarrow D_g \}_{g\in G}$ if for all $g\in G$, $D_g$ is an ideal of $\A$ and $\alpha_g$ is an isomorphism of $k$-algebras such that for all $g,h \in G$:
\begin{enumerate}
 \item[(i)] $D_e = \A$ and $\alpha_e$ is the identity map of $\A$;
\item[(ii)] $\alpha_g(D_{g^{-1}} \cap D_h) = D_g \cap D_{gh}$;
\item[(iii)] $\alpha_g(\alpha_h(x)) = \alpha_{gh}(x)$ for all $x\in D_{h^{-1}}\cap D_{(gh)^{-1}}$.
\end{enumerate}

The partial skew group ring of $\A$ and $G$ is defined to be the projective left $\A$-module $\AG = \bigoplus_{g\in G} D_g$ with multiplication
$$ (a{\o g}) (b{\o h}) = \alpha_g(\alpha_{g^{-1}}(a)b){\o {gh}}$$
for all $a\in D_g$ and $b\in D_h$ and where $\o{g}$ is the placeholder for the $g$th component of $\bigoplus_{g\in G} D_g$.  Since $\AG$ is naturally $G$-graded, the question arises how much of the Cohen-Montgomery duality carries over to partial group actions. 

As in \cite{DokuchaevFerreroPaques} we will assume that the ideals $D_g$ are generated by central idempotents, i.e. $D_g=\A1_g$ with central idempotent $1_g\in \A$ for all $g\in G$. For any $g\in G$ we define the following endomorphism $\beta_g:A\rightarrow A$ of $\A$ by $$\beta_g (a) = \alpha_g ( a 1_{g^{-1}} ) \:\:\:\forall a \in \A$$
This map gives rise to a $k$-linear map $k[G] \otimes A \rightarrow A$ with 
$$ g \otimes a \mapsto g\cdot a := \beta_g(a) = \alpha_g ( a 1_{g^{-1}} )$$
for all $g\in G, a \in \A$.

\begin{lem}\label{properties}
 With the notation above we have that
\begin{enumerate}
 \item $\beta_g$ are $k$-algebra endomorphisms of $\A$ for all $g\in G$, i.e.
$$ g \cdot (ab) = (g\cdot a)(g\cdot b) \:\:\: \forall a,b \in \A.$$
\item $g\cdot (h \cdot a) = ((gh)\cdot a) 1_g$ for all $g,h \in G$ and $a\in \A$.
\item $(g\cdot a)b = g\cdot (a (g^{-1} \cdot b))$ for all $a,b \in \A$ and $g\in G$.
\end{enumerate}
\end{lem}

\begin{proof}
 (1) follows since the $\alpha_g$ are algebra homomorphisms and the idempotents $1_g$ are central, i.e. for all $a,b \in \A$:
$$\beta_g(ab) = \alpha_g(ab1_{g^{-1}}) = \alpha_g(a1_{g^{-1}} b1_{g^{-1}}) = \alpha_g(a1_{g^{-1}})\alpha_g(b1_{g^{-1}}) = \beta_g(a)\beta_g(b). $$

(2) follows from \cite[2.1(ii)]{DokuchaevFerreroPaques}:
$$\alpha_g(\alpha_h(a1_{h^{-1}})1_{g^{-1}}) = \alpha_{gh}(a1_{h^{-1}g^{-1}}) 1_g$$
what expressed by $\beta$ yields the statement of (2).

(3) Using (1), (2) and the fact that $\beta_e=id$ and that the image of $\beta_g$ is $D_g=A1_g$ we have that
$$g\cdot (a (g^{-1} \cdot b)) = (g\cdot a)(g\cdot (g^{-1} \cdot b)) = (g\cdot a)b1_g = (g\cdot a)b.$$
\end{proof}

Obviously we also have $g\cdot 1 = \alpha_g(1_{g^{-1}})=1_g$ and $g\cdot (g^{-1} \cdot a) = ((gg^{-1})\cdot a)1_g = a1_g$
for all $a\in \A$ and $g\in G$ using property (2). Moreover using the fact that $\alpha_g$ is bijective and $1_g$ central we have for all $a\in \A$ and $g \in G$ that $g\cdot a = 0$ if and only if $a\in \A(1-1_g)$.

\section{Grading of the partial skew group ring}

The partial skew group ring  is the projective left $\A$-module $\AG = \bigoplus_{g\in G} D_g$. We will write an element of $\AG$ as a finite sum of elements $\sum_{g\in G} a_g{\o g}$ where $a_g \in D_g=A1_g$ and ${\o g}$ is a placeholder for the $g$-th component. $\AG$ becomes an associative $k$-algebra by the product:
$$ (a{\o g}) (b{\o h}) = \alpha_g(\alpha_{g^{-1}}(a)b){\o {gh}}$$ for all $g,h \in G$ and $a\in D_g$ and $b\in D_h$.  Using our '$\cdot$'-notation we see easily $$(a{\o g})(b{\o h}) =  a (g\cdot b){\o {gh}}.$$ The algebra $\AG$ is naturally $G$-graded where the homogeneous elements are those in $\{ D_g \}_{g\in G}$, i.e. $D_g D_h \subseteq D_{gh}$ by definition of the multiplication in $\AG$. Thus $\AG$ becomes a $k[G]$-comodule algebra. Note that the $G$-grading is strong, in the sense that $D_gD_h = D_{gh}$ if and only if $D_g=\A$ for all $g\in G$, i.e. the $G$-action is global (since if $D_gD_h=D_{gh}$ for all $g,h \in G$, then 
$$\A1_g1_{g^{-1}} = D_gD_{g^{-1}} = D_{gg^{-1}} = D_e = \A,$$
thus $1_g$ is an invertible central idempotent and hence equals $1$, i.e. $D_g=\A$). Known results on graded rings can be applied to the $G$-grading of $\AG$. 

\section{Duality for partial actions of finite groups}

Assume $G$ to be finite, then $k[G]^*$ becomes a Hopf algebra with projective basis $p_g \in k[G]^*$ where $p_g(h)=\delta_{g,h}$ for all $g, h\in H$. The multiplication is defined as $p_g * p_h = \delta_{g,h} p_g$ and the identity element of $k[G]^*$ is $1 = \sum_{h\in H} p_h$. Now $\AG$ becomes a $k[G]^*$-module algebra by $$p_h\triangleright  (a{\o g}) = \delta_{g,h} a{\o g} $$ for all $g,h \in G$ and $a_g\in D_g$. The multiplication of the smash product $\AGG$  is defined as
$$ (a{\o g} \# p_h) (b{\o k} \# p_l) = \sum_{s\in G} (a{\o g}) [p_{s}\triangleright (b{\o k})] \# p_{s^{-1}h}*p_l = (a{\o g})(b{\o k}) \# p_{k^{-1}h}*p_l = a(g\cdot b){\o {gk}} \# \delta_{h,kl} p_l.$$
The identity element of $\B=\AG\#k[G]^*$ is $\sum_{h\in G} 1{\o e} \# p_h$. In the case of global actions Cohen and Montgomery proved in \cite{CohenMontgomery} that $\A*G\#k[G]^* \simeq M_n(\A)$ where $n=|G|$ and $M_n(\A)$ denotes the ring of $n\times n$-matrizes over $\A$. We will index the matrizes of $M_n(\A)$ by elements of $G$ and denote by $E_{g,h}$ the elementary matrix that has the value $1$ in the $g$-th row and the $h$-th column and zero elsewhere.

\begin{prop} Let $G$ be a finite group of $n$ elements, acting partially on an $k$-algebra $\A$ and consider the $k$-algebra $\B=\AGG$. The map $$\Phi:\B \longrightarrow M_n(\A) \:\:\mbox{ with }$$ 
$$ \sum_{g,h} a_{g,h} {\o g}\# p_h \mapsto \sum_{g,h} h^{-1}\cdot ( g^{-1} \cdot a_{g,h}) E_{gh,h}$$
is a $k$-algebra homomorphism.
\end{prop}

\begin{proof}
First note that for any $g,h,k \in G$ and $a \in D_g, b \in D_h$ we have, using Lemma \ref{properties}(2) in the 2nd, 4th and 6th line and Lemma \ref{properties}(1) in the 3rd line:
\begin{eqnarray*}
 k^{-1}\cdot((gh)^{-1}\cdot (a(g\cdot b))) &=& k^{-1}\cdot \left( ((gh)^{-1}\cdot a)((gh)^{-1}\cdot (g\cdot b))\right) \\
&=& k^{-1}\cdot \left[ ((gh)^{-1}\cdot a) (h^{-1}\cdot b) 1_{(gh)^{-1}} \right] \\
&=& \left[ k^{-1}\cdot ((gh)^{-1}\cdot a)\right] \left[ k^{-1} \cdot (h^{-1}\cdot b)\right] \\
&=& ((ghk)^{-1}\cdot a) ((hk)^{-1} \cdot b) 1_{k^{-1}} \\
&=& ((ghk)^{-1}\cdot a) 1_{(hk)^{-1}} ((hk)^{-1} \cdot b) 1_{k^{-1}}  \\
&=& ((hk)^{-1}\cdot (g^{-1}\cdot a)) (k^{-1} \cdot (h^{-1} \cdot b)) 
\end{eqnarray*}
Thus we showed:
\begin{equation}\label{equation1}
 k^{-1}\cdot((gh)^{-1}\cdot (a(g\cdot b))) = ((hk)^{-1}\cdot (g^{-1}\cdot a)) (k^{-1} \cdot (h^{-1} \cdot b)) 
\end{equation}

For any $a{\o g}\#p_h,  b{\o k}\#p_l \in \AGG$ we have, using equation (\ref{equation1}):
\begin{eqnarray*}
 \Phi((a{\o g}\#p_h)(b{\o k}\#p_l)) &=& \Phi(a(g\cdot b){\o {gk}}  \# \delta_{h,kl} p_l)\\
&=&  l^{-1}\cdot((gk)^{-1}\cdot (a(g\cdot b)))E_{gkl,l} \delta_{h,kl} \\
&=& ((kl)^{-1}\cdot (g^{-1}\cdot a)) (l^{-1} \cdot (k^{-1} \cdot b)) E_{gh,h}E_{kl,l} \delta_{h,kl}\\
&=& (h^{-1}\cdot (g^{-1}\cdot a))E_{gh,h} \: (l^{-1} \cdot (k^{-1} \cdot b))E_{kl,l}\\
&=& \Phi(a{\o g}\#p_h)\Phi(b{\o k}\#p_l)
\end{eqnarray*}
Hence $\Phi$ is an algebra homomorphism.
\end{proof}

 Note that $\Phi$ restricted to $\AG$ is injective, i.e. $\AG$ can be considered a subalgebra of $M_n(\A)$.
In general $\Ker{\Phi}$ is non-trivial, unless the partial action is a global action.

\begin{prop}
 $\Ker{\Phi} = \bigoplus_{g,h \in G} \A(1-1_{gh})1_g {\o g} \# p_h$.
\end{prop}
\begin{proof}
 Suppose $\gamma =\sum_{g,h} a_{g,h}{\o g} \# p_h \in \Ker{\Phi}$, then $h^{-1}\cdot (g^{-1} \cdot a_{g,h})=0$ for all $g,h \in G$. Thus $(g^{-1} \cdot a_{g,h})\in \A(1-1_h)\cap D_{g^{-1}}=\A(1-1_h)1_{g^{-1}}$. Hence 
$$ a_{g,h} = g \cdot (g^{-1} \cdot a_{g,h}) \in \A g\cdot(1 - 1_h) = \A (1_g - 1_g1_{gh}), $$
i.e. $\gamma \in \bigoplus_{g,h} \A(1-1_{gh})1_g{\o g}\# p_h$. The other inclusion follows because
$$\Phi\left((g\cdot (1-1_h)){\o g}\# p_h\right) = h^{-1}\cdot (g^{-1}\cdot (g\cdot (1-1_h))) E_{gh,h} =  h^{-1}\cdot ((1-1_h)1_g) E_{gh,h}  = 0.$$
\end{proof}

Note that the inclusion of $\AG$ into $\AGG$ is given by $a {\o g} \mapsto \sum_{h\in G} a {\o g} \# p_h$ for all $g\in G$ and $a \in D_g$. If $\sum_{h\in G} a{\o g} \# p_h \in \Ker{\Phi}$, then $a \in \A(1-1_{gh})1_g$ for all $h\in G$. In particular for $h=e$ we have $a \in \A(1-1_g)1_g = 0$. Hence $\Phi$ restricted to $\AG$ is injective.

We will describe the image of $\Phi$. By definition of $\Phi$, the image of an arbitrary element $\gamma=\sum_{g,h} a_{g,h}{\o g} \# p_h$ is 
$$ \Phi(\gamma) = \sum_{g,h} ( (gh)^{-1} \cdot a_{g,h}) 1_{(gh)^{-1}} 1_{h^{-1}} E_{gh,h} = {\left(b_{r,s}1_{r^{-1}}1_{s^{-1}}\right)}_{r,s\in G}$$
with $b_{r,s} =  r^{-1} \cdot a_{rs^{-1}, s} $ for all $r,s\in G$.

\begin{prop}
The image of $\Phi$ consists of all matrices of the form 
${\left(b_{g,h}1_{g^{-1}}1_{h^{-1}}\right)}_{g,h\in G}$ for any matriz $(b_{g,h})$ of elements of $\A$. In particular
$\Im{\Phi} = \e M_n(A) \e$, where $\e$ is the idempotent $\sum_{g\in G} 1_{g^{-1}} E_{g,g}$.
\end{prop}
\begin{proof}
 We saw already that an element of the image of $\Phi$ is of the given form. 
Note that by definition of partial group action we have 
$$ D_g \cap D_{gh} = \alpha_g(D_{g^{-1}} \cap D_h)$$
for all $g,h \in G$. Hence also 
$$ D_{g^{-1}} \cap D_{h^{-1}} = \alpha_{g^{-1}}(D_g \cap D_{gh^{-1}})$$
holds for all $g,h \in G$. Thus for all $ b\in \A$ there exists $a\in \A$ such that 
$$b1_{g^{-1}}1_{h^{-1}} = \alpha_{g^{-1}} (a 1_{gh^{-1}} 1_g) = g^{-1} \cdot (a1_{gh^{-1}}).$$
This implies that
\begin{eqnarray*}
\Phi(a1_g1_{gh^{-1}} {\o {gh^{-1}}} \# p_h) &=&
h^{-1}\cdot ((hg^{-1})\cdot  (a1_g1_{gh^{-1}})) E_{g,h} \\
&=& g^{-1} \cdot  (a1_g1_{gh^{-1}})) 1_{h^{-1}} E_{g,h} \\
&=& b 1_{g^{-1}}1_{h^{-1}} E_{g,h} 
\end{eqnarray*}
Hence given any matrix $(b_{g,h})$ there are elements $a_{g,h}$ such that 
$$\Phi\left(\sum_{g,h} a_{g,h}1_g1_{gh^{-1}} {\o {gh^{-1}}} \# p_h\right) =  \sum_{g,h} b_{g,h} 1_{g^{-1}}1_{h^{-1}} E_{g,h} =  {\left(b_{g,h}1_{g^{-1}}1_{h^{-1}}\right)}_{g,h\in G}.$$
This shows that $\Im{\Phi}$ consists of all matrizes of the given form and hence is equal to $\e M_n(A) \e$. Note that $\e$ is the image of the identity element of $\B$.
\end{proof}

The last Propositions yield our main result in this section

\begin{thm}\label{thetheorem}
 $\AGG \simeq \Ker{\Phi} \times \e M_n(\A) \e$.
\end{thm}

\begin{proof}
 The kernel of $\Phi$ is an ideal and a direct summand of $\B=\AGG$. To see this we first show that the left $\A$-module $I=\bigoplus_{g,h \in G} \A1_{gh}1_g{\o g} \# p_h$ is a two-sided ideal of $\B$. 
For any $x{\o k}\#p_l \in \B$ and $a1_{gh}1_g {\o g} \# p_h \in I$ we have 
\begin{eqnarray*}
(a1_{gh}1_g {\o g} \# p_h) (b{\o k}\#p_l) &=& a1_{gh} 1_g (g \cdot b1_k){\o {gk}} \# \delta_{h,kl} p_l = a(g \cdot b)\delta_{h,kl} 1_{gkl} 1_{gk} {\o {gk}} \# p_l \in I.\\
(b{\o k}\#p_l)(a1_{gh}1_g {\o g} \# p_h) &=& b (k\cdot a1_{gh} 1_g) {\o {kg}} \# \delta_{k,gh}p_h = b(g \cdot a)\delta_{h,kl} 1_{kgh} 1_{kg} {\o {kg}} \# p_h \in I.
\end{eqnarray*}
Since $I\oplus \Ker{\Phi} = \B$ and both direct summands are two-sided ideals we have $\B = I\times \Ker{\Phi}$ (ring direct product). Moreover $\Phi(I)=\e M_n(\A) \e = \Im{\Phi}$. This implies $\B\simeq \Ker{\Phi} \times \e M_n(\A) \e$.

\end{proof}

Note that $\Phi$ embedds $\AG$ into the Pierce corner $\e M_n(\A) \e$.
\begin{cor}
 $\AG$ is isomorphic to a separable subalgebra of $\e M_n(\A) \e$. 
\end{cor}
\begin{proof}
Recall that the subalgebra $\AG$ sits into $\B$ by $a{\o g} \mapsto \sum_{h\in G} a{\o g} \#p_h$.
The right action of $\AG$ on $\B$ is given by 
$$(x{\o k} \#p_l)\cdot a{\o g} =  (x{\o k} \#p_l)\left(\sum_{h\in G} a{\o g} \#p_h\right) = (x{\o k})(a{\o g})\#  p_{g^{-1}l}$$
The left action is given by
$$a{\o g} \cdot (x{\o k} \#p_l) = \left(\sum_{h\in G} a{\o g} \#p_h\right) (x{\o k} \#p_l) = (a{\o g})(x{\o k})\#  p_l$$
The element
$$f=\sum_{g\in G} {\o e}\# p_g \otimes {\o e} \# p_g \in \B \otimes_{\AG} \B$$ is $\AG$-centralising, i.e.
for all $a{\o h} \in \AG$ we have
$$ f a{\o h} = \sum_{g\in G} {\o e}\#p_g \otimes a{\o h} \#p_{h^{-1}g} = \sum_{g\in G} a{\o h} \#p_{h^{-1}g} \otimes {\o e} \#p_{h^{-1}g} = a{\o h} f$$
Since also $\mu(f)= {\o e} \# \sum_{g\in G} p_g = 1_{\B}$ we have that $f$ is a seperability idempotent for $\B$ over $\AG$.
Hence $\e M_n(\A) \e \simeq  \Phi(\B)$ is separable over $\Phi(\AG) \simeq \AG$.
\end{proof}

\section{Trivial partial actions}
The easiest example of partial actions arise from (central) idempotents in a $k$-algebra $\A$. Suppose that $\A$ admits a non-zero central idempotent, i.e. there exist algebras $R,S$ such that $\A=R\times S$ as algebras. For any group $G$ set $D_g = R\times 0$ and $\alpha_g = id_{D_g}$for all $g\neq e$ and $D_e=\A$ and $\alpha_e=id_{\A}$.  Then $\{ \alpha_g \mid g \in G\}$ is a partial action of $G$ on $\A$. The partial skew group ring turns out to be $\AG \simeq R[G] \times S$, where $R[G]$ denotes the group ring of $R$ and $G$. Note that $0\times S$ is in the zero-componente of the $G$-grading on $\AG$. If $G$ is finite, say of order $n$, then a short calculation (using Cohen-Montgomery duality and Theorem \ref{thetheorem}) shows that $\B = (\AG)\#k[G]$ is isomorphic to $M_n(R) \times S^n$ where $S^n$ denotes the direct product of $n$ copies of $S$. Depending on the rings $R$ and $S$, $\B$ might or might not be Morita equivalent to $\A$. For instance if $R=S=F$ is a field, then any progenerator $P$ for $\A$ has the form $F^k \times F^m$ for numbers $k,m \geq 1$. Thus $\End{P} \simeq M_k(F) \times M_m(F)$, whose center is isomorphic to $F^2 = \A$. On the other hand $\B=(\AG)\#k[G] \simeq M_n(F)\times F^n$ has center $F^{n+1}$, i.e. $\B$ will be Morita equivalent to $\A$ if and only if $G$ is trivial.

On the other hand, there are algebras which satisfy (as algebras) $\A^n \simeq \A \simeq M_n(\A)$ for any $n$. To give an example, let $R$ be the ring of sequences of elements of a field $k$, i.e. $R=k^{\NN}$. The function $\chi$ with $\chi(2n)=1$ and $\chi(2n+1)=0$ for all $n$  defines an idempotent of $R$. The map $\Psi: \chi R \rightarrow R$ with $\Psi(\chi f)(n)=f(2n)$ is a ring isomorphism. Analogosuly we can show that $(1-\chi)R \simeq R$. Hence $R^2\simeq R$.
Now take $\A = \End{F}$, where $F=R^{(\NN)}$ denotes the countable infinite free $R$-module. Using again $\chi$ we have that $$\A = (\chi \A) \times ((1-\chi) \A) \simeq \A \times \A \simeq \cdots \simeq  \A^n$$
for any $n\geq 2$. Moreover for any partition of $\NN$ into $n$ infinite disjoint subsets $\Lambda_1, \ldots, \Lambda_n$, we have that $$F=R^{(\NN)} \simeq R^{(\Lambda_1)} \oplus \cdots \oplus R^{(\Lambda_n)} \simeq F^n.$$
Hence $\A = \End{F} \simeq \End{F^n} \simeq M_n(\A)$. Applying the double skew group ring construction
 again we conclude that 
$$ \B = (\AG)\# k[G] \simeq M_n(\chi\A) \times ((1-\chi)\A)^n \simeq \A \times \A \simeq \A.$$


\section{Infinite partial group action}

Following Quinn \cite{quinn} we define $\Phi$ in case of $G$ being infinite as a map from $\AG$ to the ring of row and column finite matrizes. Let $M_G(\A)$ be the subring of $\End{\A^{|G|}}$ consisting of row and column finite matrizes $(a_{g,h})_{g,h \in G}$ indexed by elements of $G$ with entries in $\A$, i.e. for any $g\in G$ the sets $\{ a_{gh} | h\in G\}$ and $\{a_{hg} | h \in G\}$ are finite. Let $E_{g,h}$ be, as above, those matrizes that are $1$ in the $(g,h)$th component and zero elsewhere. Note that $E_{g,h}E_{r,s}=\delta_{h,r}E_{g,s}$. Then define $\Phi: \AG \rightarrow M_G(\A)$ by 
$$ a{\o g} \mapsto \sum_{h \in G} h^{-1}\cdot (g^{-1}\cdot a) E_{gh,h}$$
for any $a{\o g} \in \AG$. Note that the (infinite) sum on the right side makes sense in $M_G(\A)$. As above one checks that $\Phi$ is an algebra homomorphism.

\begin{prop}
 Let $G$ be any group acting partially on $\A$. Then $\AG$ is isomorphic to a subalgebra of $\e M_G(A) \e$ where $M_G(\A)$ denotes the ring of row and column finite matrizes indexed by elements of $G$ and with entries in $\A$. The element $\e$ is the idempotent $\sum_{g\in G} 1_{g^{-1}} E_{g,g}$.
\end{prop}

\begin{proof}
For all $a{\o g}, b{\o h} \in \AG$ we have using equation (\ref{equation1}) in the 4th line:
\begin{eqnarray*}
 \Phi(a{\o g})\Phi(b{\o h}) &=& \left(\sum_{k\in G} k^{-1}\cdot (g^{-1} \cdot a) E_{gk,k} \right) \left( \sum_{l\in G} l^{-1}\cdot (h^{-1} \cdot b) E_{hl,l}\right)\\
&=& \sum_{k,l \in G} (k^{-1}\cdot (g^{-1} \cdot a)) (l^{-1}\cdot (h^{-1} \cdot b)) E_{gk,k} E_{hl,l}\\
&=& \sum_{l \in G} ((hl)^{-1}\cdot (g^{-1} \cdot a)) (l^{-1}\cdot (h^{-1} \cdot b)) E_{ghl,l}\\
&=& \sum_{l \in G} l^{-1} \cdot ((gh)^{-1}\cdot (a (g\cdot b))) E_{ghl,l}\\
&=& \Phi(a (g\cdot b){\o gh}) \\
&=& \Phi((a{\o g})(b{\o h})) 
\end{eqnarray*}
Hence $\Phi$ is an algebra homomorphism. 
Since 
$$\Phi(a{\o g})=0 \Leftrightarrow (\forall h\in G): h^{-1}\cdot (g^{-1} \cdot a) = 0 \Rightarrow g\cdot (g^{-1}\cdot a) = a1_g = 0 \Rightarrow a=0,$$ we have that $\Phi$ is injective. Moreover $\Phi(a{\o g})\in \e M_G(A) \e$ as above.
\end{proof}

\section{Partial Hopf action}

In \cite{CaenepeelJanssen} Caenepeel and Janssen defined the notion of a partial Hopf action as follows:
Let $H$ be a Hopf algebra, with comultiplication $\Delta$, counit $\epsilon$ and antipode $S$, and let $\A$ be a $k$-algebra such that there exists a $k$-linear map 
$$\cdot: H\otimes A \rightarrow A$$ sending $h\otimes a \mapsto h\cdot a$.
The action $\cdot$ is called a {\it partial Hopf action} if for all $h,g \in H$ and $a,b \in \A$:
\begin{enumerate}
 \item[(1)] $h \cdot (ab) = \sum_{(h)} (h_1\cdot a)(h_2 \cdot b)$;
\item[(2)] $1_H \cdot a = a $;
\item[(3)] $h\cdot (g\cdot a) = \sum_{(h)} (h_1 \cdot 1)((h_2g)\cdot a)$;
\end{enumerate}

Let $H$ be a Hopf algebra which is finitely generated and projective as $k$-module with dual basis $\{ (b_i,p_i)\in H\times H^* \mid 1\leq i\leq n\}$. Then there exist structure constants $c_{k,l}^i$ and $m_{k,l}^i$ in $k$ such that
$\Delta(b_i)=\sum_{k,l=1}^n c_{k,l}^i b_k\otimes b_l$ and
$b_kb_l = \sum_{i=1}^n m_{k,l}^i b_i$ for all $1\leq i,k,l \leq n$.
It is well-known that $H^*$ becomes a Hopf algebra with comultiplication and multiplication defined on the generators $\{p_i \mid 1\leq i \leq n\}$ as follows:
$\Delta_{H^*}(p_i) = \sum_{k,l=1}^n m_{k,l}^i p_k \otimes p_l$
and $p_k*p_l = \sum_{i=1}^n c_{k,l}^i p_i$. The counit of $H^*$ is given by $\epsilon_{H^*}(f) = f(1)$.

Recall that $H^*$ acts on $H$ from the left by $f \rhu h = \sum_{(h)} h_1 f(h_2)$, such that the smash product $H\# H^*$ can be considered whose multiplication is given by $$ (h\# f)(k\# g)= \sum_{(f)} h(f_1 \rhu k) \# f_2 * g$$
for all $h,k \in H$ and $f,g \in H^*$. The smash product yields a left module action on $H$, i.e. an algebra homomorphism
$$ \lambda: H \# H^* \rightarrow \End{H} \:\:\: h\# f \mapsto [k \mapsto h(f\rhu k)].$$

The smash product $H\# H^*$ is sometimes called the Heisenberg double of $H$ and in case $H$ is free of finite rank isomorphic to $\End{H}$ (see \cite[9.4.3]{montgomery}).

Analougosly we have a right action of $H^*$ on $H$ by $h\lhu f = \sum_{(h)} h_2 f(h_1)$ for all $f\in H^*$ and $h\in H$, turning $H$ into a right $H^*$-module algebra. The smash product $H^*\# H$ yields a right module action on $H$, i.e. an algebra anti-homomorphism
$$ \rho: H^* \# H \rightarrow \End{H} \:\:\: f\# h \mapsto [k \mapsto (k\lhu f)h]$$
As in \cite[9.4.10]{montgomery} one shows that for all $h,k\in H$ and $f,g \in H^*$:
\begin{equation}\label{equation_montgomery}
 \lambda(h\# f)\rho(g\# 1) = \sum_{(g)} \rho(g_2 \# 1)\lambda((h\lhu S(g_1))\# f) 
\end{equation}
Now assume that $H$ acts partially on $\A$, then the map
$\Delta_{\A} : A \rightarrow A \otimes H^*$ with
$$ \Delta(a)_A =  \sum_{i=1}^n (b_i\cdot a) \otimes p_i$$
for all $a\in \A$ defines a {\it partially coaction}. The map $\Delta_{\A}$ satisfies:
\begin{eqnarray*}
\Delta_{\A}(ab) &=& \Delta_{\A}(a) \Delta_{\A}(b)\\
(1\otimes \epsilon_{H^*})\Delta_{\A}(a) &=& id_A(a)\\
(\Delta_{\A}\otimes 1)\Delta_{\A}(a) &=& (\Delta_{\A}(1)\otimes 1)(1 \otimes \Delta_{H^*})\Delta_{\A}(a)
\end{eqnarray*}
The last equation shows that in general this coaction does not make $\A$ into a right $H$-comodule. It  can be deduced using the structre constants and property (3) from above
\begin{eqnarray*}
 \sum_{i,j=1}^n b_j \cdot (b_i\cdot a) \otimes p_j \otimes p_i &= &
 \sum_{i,j,k,l,r=1}^n c_{k,l}^j m_{l,i}^r (b_k\cdot 1)(b_r \cdot a) \otimes p_j \otimes p_i\\
&= &  \sum_{i,k,l,r=1}^n m_{l,i}^r (b_k\cdot 1)(b_r \cdot a) \otimes p_kp_l \otimes p_i\\
&= &  \sum_{k,r=1}^n (b_k\cdot 1)(b_r \cdot a) \otimes (p_k \otimes 1)\Delta_{H^*}(p_r)\\
&=& \left(\sum_{k=1}^n (b_k \cdot 1) \otimes p_k\right) \left(\sum_{r=1}^n (b_r\cdot a)\otimes \Delta(p_r)\right)
\end{eqnarray*}
With the above notation we define a homomorphism $\phi:\A \rightarrow \A\otimes \End{H}$ by
$$\phi(a) = \sum_{i=1}^n (b_i\cdot a) \otimes \rho(S^{-1}(p_i) \# 1).$$
Then $\phi$ is an algebra homomorphism, because
\begin{eqnarray*}
\phi(ab)&=&\sum_{i=1}^n (b_i\cdot (ab) \otimes \rho(S^{-1}(p_i) \# 1)\\
&=&\sum_{i=1}^n ((b_i)_1\cdot a)(((b_i)_2\cdot b) \otimes \rho(S^{-1}(p_i) \# 1)\\
&=&\sum_{k,l=1}^n (b_k\cdot a)((b_l\cdot b) \otimes \rho(S^{-1}(c_{k,l}^ip_i) \# 1)\\
&=&\sum_{k,l=1}^n (b_k\cdot a)((b_l\cdot b) \otimes \rho(S^{-1}(p_l)S^{-1}(p_k) \# 1)\\
&=&\sum_{k,l=1}^n (b_k\cdot a)((b_l\cdot b) \otimes \rho(S^{-1}(p_k) \# 1) \rho(S^{-1}(p_l)\# 1) \\
&=& \phi(a)\phi(b).
\end{eqnarray*}
where we use in the line before the last the fact that $\rho$ is an anti-homomorphism.

The partial smash product of $\A$ and $H$ is defined as a certain submodule of $\A\otimes H$.
On $\A\otimes H$ we define a new (associative) multiplication by 
$$ (a\otimes h)(b\otimes g) := \sum_{(h)} a (h_1\cdot b) \otimes h_2g.$$
for all $a,b \in \A$, $h,g \in H$.
Note that $\A\otimes H$ is naturally an $\A$-bimodule given by 
$$x (a\otimes h) y = (x\otimes 1)(a\otimes h)(y\otimes 1) = \sum_{(h)} xa(h_1\cdot y) \otimes h_2$$
The partial smash product is defined to be $\AH = (\A\otimes H) 1_{\A}$ and is spanned by the elements of the form $\sum_{(h)} a (h_1\cdot 1_{\A}) \otimes h_2 \}$ for all $a\in \A, h \in H$.
The partial smash product becomes naturally a right $H$-comodule algebra by $$\rho = 1\otimes \Delta : \A\otimes H \rightarrow \A\otimes H \otimes H, \:\:\: a\otimes h \mapsto \sum_{(h)} a \otimes h_1 \otimes h_2$$
and for all $(a\otimes h)1_{\A} \in \AH$ we have $$ \rho( (a\otimes h)1_{\A} ) =  \sum_{(h)} a (h_1 \cdot 1_{\A}) \otimes h_2 \otimes h_3,$$  making $\AH$ into a right $H$-comodule algebra. Moreover $\AH$ becomes a left $H^*$-module algebra, where the action is defined by $$f \rhd ( (a\#h)1_{\A}) = \sum_{(h)} (a(h_1\cdot 1_{\A}) \# (f\rhu h_2) = (a \# (f\rhu h))1_{\A},$$ for all $f\in H^*, h \in H, a\in \A$.
The classical Blattner-Montgomery duality (\cite{BlattnerMontgomery} says that the double smash product $\A\#H\#H^*$ is isomorphic to $M_n(\A)$ where $n$ is the rank of $H$ over $k$.

\begin{lem}
Let $\psi:H\# H^* \rightarrow \A\otimes \End{H}$ be the map defined by $h\# f \mapsto 1 \otimes \lambda(h\# f)$ for all $h\in H, f \in H^*$. Then for all $a\in \A, h \in H, f \in H^*$ we have
$$ \phi(1)\psi(h\# f)\phi(a) = \sum_{(h)} \phi(h_1\cdot a)\psi(h_2 \# f).$$
\end{lem}

\begin{proof}
Let $a\in \A, h \in H, f \in H^*$.
 \begin{eqnarray*}
\sum_{(h)} \phi(h_1\cdot a)\psi(h_2 \# f) &=& \sum_{(h),i} p_i(h_1) \phi(b_i\cdot a)\psi(h_2 \# f)\\
&=& \sum_{i,j} b_j\cdot (b_i\cdot a) \otimes \rho(S^{-1}(p_j)\#1)\lambda(h\lhu p_i  \# f)\\
&=& \sum_{k,r} (b_k\cdot 1)(b_r\cdot a) \otimes \rho(S^{-1}(p_k(p_r)_1)\#1)\lambda(h\lhu (p_r)_2  \# f)\\
&=& \sum_{k,r} (b_k\cdot 1)(b_r\cdot a) \otimes \rho(S^{-1}(p_k) \# 1) \rho(S^{-1}((p_r)_1)\#1)  \lambda(h\lhu (p_r)_2  \# f)\\
&=& \phi(1) \sum_{r} (b_r\cdot a) \otimes \rho(S^{-1}(p_r)_2 \#1)\lambda(h\lhu (S(S^{-1}(p_r)_1)  \# f)\\
&=& \phi(1) \sum_{r} (b_r\cdot a) \otimes \lambda(h\# f) \rho(S^{-1}(p_r) \# 1)\\
&=& \phi(1) \psi(h\# f)\phi(a)
 \end{eqnarray*}
where we  use equation (\ref{equation_montgomery}) in the third line from below.
\end{proof}

\begin{thm}
 Suppose that $H$ is a Hopf algebra, finitely generated projective over $k$, which partially actions on $\A$. Then $\Phi:  A \otimes H \# H^* \rightarrow \A\otimes \End{H}$ with 
$$a\otimes h \# f \mapsto \phi(a)\psi(h\# f)$$ is an algebra homomorphism. The image of the restriction to $\AH \# H^*$ lies inside $\e \left( A \otimes \End{H}\right) \e$ where $\e$ is the idempotent $$\e=\sum_{i=1}^n (b_i \cdot 1) \otimes \rho(S^{-1}(p_i)\otimes 1).$$
\end{thm}

\begin{proof}
 For any $a,b \in \A, h,k \in H$ and $f,g \in H^*$ we have 
\begin{eqnarray*}
\Phi(a\otimes h \# f) \Phi(b\otimes k \# g) ) &=& \phi(a)\psi(h\# f)\phi(b)\psi(k\# g)\\
&=& \phi(a)\phi(1)\psi(h\# f)\phi(b)\psi(k\# g)\\
&=& \sum_{(h)}\phi(a)\phi(h_1 \cdot b)\psi(h_2\# f)\psi(k\# g)\\
&=& \sum_{(h,f)}\phi(a(h_1 \cdot b))\psi(h_2(f_1\rhu k) \# f_2*g)\\
&=& \Phi\left(\sum_{(h,f)} a(h_1 \cdot b) \otimes h_2(f_1\rhu k) \# f_2*g\right)\\
&=& \Phi\left((a\otimes h \# f)(b\otimes k \# g)\right).
\end{eqnarray*}
Hence $\Phi$ is an algebra homomorphism. Since the image of the identity $\mathbf{1}=1_{\A}\#1_H\#1_{H^*}$ of $\AH\#H^*$ under the map $\Phi$ is $\e$, $\e$ is an idempotent. Moreover 
$$\Phi(\gamma) = \Phi ( \mathbf{1} \gamma \mathbf{1}) \in \e (\A\otimes \End{H})\e, $$
for all $\gamma \in \AH\# H^*$.
\end{proof}


\begin{thebibliography}{CWZ}

\bibitem{BlattnerMontgomery} R.J.Blattner and S.Montgomery, \textit{A duality theorem for Hopf module algebras.}
J. Algebra 95 (1985), no. 1, 153--172. 


\bibitem{CaenepeelJanssen} S.Caenepeel and  K. Janssen \textit{Partial entwining structures}, arxiv.org/math/0610524

\bibitem{CohenMontgomery} M.Cohen and S.Montgomery, \textit{Group-graded rings, smash products, and group actions.},
Trans. Amer. Math. Soc. 282 (1984), no. 1, 237--258. 

\bibitem{CortesFerrero} W.Cortes and M.Ferrero, \textit{Partial skew polynomial rings: prime and maximal ideals.}, Comm. Algebra  35  (2007),  no. 4, 1183--1199.
  
\bibitem{DokuchaevFerreroPaques} M. Dokuchaev, M.Ferrero and A.Paques, \textit{Partial actions and Galois theory.}, J. Pure Appl. Algebra 208 (2007), no. 1, 77--87. 

\bibitem{exel} R. Exel, \textit{ Circle actions on $C^*$-algebras, partial automorphisms and generalized Pimsner-
Voiculescu exact sequences}, J. Funct. Anal. 122 (1994), 361?401.

\bibitem{montgomery} S.Montgomery, \textit{Hopf algebras and their actions on rings.}, CBMS Regional Conference Series in Mathematics, 82. AMS (1993)

\bibitem{quinn} D.Quinn,\textit{Group-graded rings and duality.}  Trans. Amer. Math. Soc.  292  (1985),  no. 1, 155--167. 


\end{thebibliography}
\end{document}